\documentclass[11pt,reqno]{amsart}
\usepackage[percent]{overpic}
\usepackage{amsfonts,amssymb,amsmath,color}

\usepackage{bbm}

\usepackage[pdftex,colorlinks,
citecolor=black,linkcolor=black,urlcolor=black]{hyperref}

\newtheorem{theorem}{Theorem}

\numberwithin{equation}{section} \numberwithin{theorem}{section}

\begin{document}

\title{Pattern Recognition on Oriented Matroids: Symmetric Cycles in the Hypercube Graphs. VI}

\author{Andrey O. Matveev}
\email{andrey.o.matveev@gmail.com}
\urladdr{https://github.com/andreyomatveev}

\begin{abstract}
We briefly discuss linear algebraic, combinatorial, and applied aspects of an exact model representation of binary arrays.
As an illustration, we present two linear algebraic portraits of a string of characters.
\end{abstract}

\maketitle

\pagestyle{myheadings}

\markboth{PATTERN RECOGNITION ON ORIENTED MATROIDS}{A.O.~MATVEEV}

\thispagestyle{empty}



\section{Introduction}

Let
\begin{equation}
\label{eqqq:2}
\boldsymbol{F}:=
\left(
\begin{smallmatrix}
T^1
\\ T^2
\\ \vdots
\\ T^{\tau}
\end{smallmatrix}
\right)
:=
\left(
\begin{smallmatrix}
(F^1)^{\top}
& (F^2)^{\top}
& \cdots
& (F^t)^{\top}
\end{smallmatrix}
\right)
\in\{1,-1\}^{\tau\times t}\subset\mathbb{R}^{\tau\times t}
\end{equation}
be a $\tau\times t$ matrix (such that $t\geq 3$)\footnote{We are mostly interested in the situation where
$\tau<\!\!< t$,
that is, the matrix~(\ref{eqqq:2}) has a very large number of columns.} whose entries
are real numbers $1$ and~$-1$, and whose rows
\begin{equation}
\label{eqqq:3}
T^i:=(T^i(1),T^i(2),\ldots,T^i(t))
\end{equation}
and columns
\begin{equation*}
(F^j)^{\top}:=(F^j(1),F^j(2),\ldots,F^j(\tau))^{\top}
\end{equation*}
are indexed starting with $1$.
We regard the rows $T^i$ of the matrix $\boldsymbol{F}$ as vectors of the real Euclidean space~$\mathbb{R}^t$. It is also convenient to interpret the row vectors~(\ref{eqqq:3}) as vertices of the {\em hypercube graph\/}~$\boldsymbol{H}(t,2)$. A pair of vertices~$\{T', T''\}\subset\{1,-1\}^t$ is an {\em edge\/} of the graph~$\boldsymbol{H}(t,2)$ if and only if the {\em Hamming distance\/}~$d(T',T''):=|\{k\in[t]\colon T'(k)\neq T''(k)\}|$ between $T'$ and~$T''$ is $1$; here $[t]$ denotes the set~$\{1,2,\ldots,t\}$.

In this note, we briefly discuss linear algebraic, combinatorial, and applied aspects of an exact model representation of binary matrices/arrays~(\ref{eqqq:2}) and of their~$\{0,1\}$-analogues. As an illustration, we present two linear algebraic portraits of a string of characters.

\section{Computation-free linear algebraic decompositions of vertices in the hypercube graph $\boldsymbol{H}(t,2)$ with respect to its  distinguished symmetric cycle}

We are interested in a distinguished {\em symmetric cycle\/} $\boldsymbol{R}:=(R^0,R^1,\ldots,$ $R^{2t-1},R^0)$ of the hypercube graph~$\boldsymbol{H}(t,2)$, with its {\em vertex sequence\/}
\begin{equation*}
\vec{\mathrm{V}}(\boldsymbol{R}):=(R^0,R^1,\ldots,R^{2t-1})
\end{equation*}
described as follows:
\begin{equation}
\label{eqqq:5}
\begin{split}
R^0:\!&=(1,1,\ldots,1)\; ,\\
R^s:\!&=(
\underbrace{-1,-1,\ldots,-1}_{s},
\underbrace{1,1,\ldots,1}_{t-s})\; ,\ \ \ 1\leq s\leq t-1\; ,
\end{split}
\end{equation}
and
\begin{equation}
\label{eqqq:6}
R^{k+t}:=-R^k\; ,\ \ \ 0\leq k\leq t-1\; .
\end{equation}

For any row~$T^i$ of the matrix~(\ref{eqqq:2}) there is a {\em unique inclusion-minimal\/} subset $\boldsymbol{Q}(T^i,\boldsymbol{R})$ of the vertex sequence~$\vec{\mathrm{V}}(\boldsymbol{R})$ such that
\begin{equation}
\label{eqqq:4}
T^i=\sum_{Q\in\boldsymbol{Q}(T^i,\boldsymbol{R})}Q\; .
\end{equation}
Since the sequence $\vec{\mathrm{V}}(\boldsymbol{R})$
is an (ordered) {\em maximal positive basis\/} of the space~$\mathbb{R}^t$, the decomposition set $\boldsymbol{Q}(T^i,\boldsymbol{R})$ of the vertex~$T^i$ of the graph~$\boldsymbol{H}(t,2)$ with respect to the cycle~$\boldsymbol{R}$ is {\em linearly independent}; its cardinality
\begin{equation*}
\mathfrak{q}(T^i):=|\boldsymbol{Q}(T^i,\boldsymbol{R})|
\end{equation*}
is {\em odd}.

Thanks to Proposition~5.9 in~\cite{M-JSP} (see~\cite[Prop.~2.4]{M-02601}), in practice, explicit decompositions of the form~(\ref{eqqq:4}) need no any computations, and thus they are fast. Moreover, since the decomposition mechanism is based solely on the inspection of the interval structure of
the {\em negative part\/} $\{e\in[t]\colon T^i(e)=-1\}$ of the vertex~$T^i$ of the graph~$\boldsymbol{H}(t,2)$, the decomposition set~$\boldsymbol{Q}(T^i,\boldsymbol{R})$ can be created even on-the-fly if the entries of the row~$T^i$ of the matrix~(\ref{eqqq:4}) themselves are revealed `sequentially'.

No instances of vectors of the maximal positive basis~$\vec{\mathrm{V}}(\boldsymbol{R})$ need to be stored into the operational memory of a `device'. The reason is that the vectors~$R^k\in\vec{\mathrm{V}}(\boldsymbol{R})$ have a simply-organized inner structure: the negative part of~$R^k$ constitutes one (empty, in the case of the vector~$R^0$) interval of the set~$[t]$. As a consequence, any component of the vector~$R^k$
can be instantly retrieved, when needed, via the mechanism of {\em lazy generation\/} based on definition~(\ref{eqqq:5})(\ref{eqqq:6}).

We see that, in practice, a report on the results of the decomposition procedure given in~(\ref{eqqq:4}) should not consist necessarily in the set~$\boldsymbol{Q}(T^i,\boldsymbol{R})$ described explicitly. Instead, if
\begin{equation}
\label{eqqq:8}
\boldsymbol{Q}(T^i,\boldsymbol{R}):=\bigl\{R^{k_0(T^i)},R^{k_1(T^i)},\ldots,R^{k_{\mathfrak{q}(T^i)-1}(T^i)}\bigr\}
\subset\vec{\mathrm{V}}(\boldsymbol{R})\; ,
\end{equation}
for some indices-superscripts $k_0(T^i)<k_1(T^i)<\cdots < k_{\mathfrak{q}(T^i)-1}(T^i)$,
then it suffices to return to the user the set of those indices
\begin{equation}
\label{eqqq:7}
\{k_0(T^i),k_1(T^i),\ldots,k_{\mathfrak{q}(T^i)-1}(T^i)\}\; .
\end{equation}

If $T^i\in\vec{\mathrm{V}}(\boldsymbol{R})$, then $\mathfrak{q}(T^i)=1$. If $T^i\not\in\vec{\mathrm{V}}(\boldsymbol{R})$, then it is shown between the lines in (the proof of)~\cite[Theorem~1.17(i)]{M-PROM} that
\begin{equation*}
\mathfrak{q}(T^i)\leq
\begin{cases}
t\; , & \text{if $t$ odd}\; ,\\
t-1\; , & \text{if $t$ even}\; .
\end{cases}
\end{equation*}
To put it simply, Proposition~5.9 in~\cite{M-JSP} (\cite[Prop.~2.4]{M-02601}) explains that the greater the number of inclusion-maximal intervals composing the negative part of~$T^i$, the larger the size of its decomposition set~$\boldsymbol{Q}(T^i,\boldsymbol{R})$, and vice versa.
Hence, we obtain the following bounds on the total number of integers (that provide us, as basic linear algebra guarantees, with a precise portrait of the matrix $\boldsymbol{F}$ defined in~(\ref{eqqq:2})) appearing in the sets~(\ref{eqqq:7}), for all of the vectors~$T^i$, $1\leq i\leq\tau$:

\begin{theorem}
\label{proppp:1}
For the tuple
\begin{equation}
\label{eqqq:12}
\begin{split}
\Big(\ \quad t\; , &\quad\ \tau\; ,\\
\{k_0(T^1),k_1(T^1),&\ldots,k_{\mathfrak{q}(T^1)-1}(T^1)\}\; ,\\
\{k_0(T^2),k_1(T^2),&\ldots,k_{\mathfrak{q}(T^2)-1}(T^2)\}\; ,\\
&\ldots\; ,\\
\{k_0(T^{\tau}),k_1(T^{\tau}),&\ldots,k_{\mathfrak{q}(T^{\tau})-1}(T^{\tau})\}\ \Big)\; ,
\end{split}
\end{equation}
consisting of the dimensionality parameters $t$ and~$\tau$,\footnote{$\quad$ Note that the parameter $\tau$ is somewhat redundant, because the tuple~(\ref{eqqq:12}) includes exactly $\tau$ sets of indices.} and of the sequence of sets of integers
described in~{\rm(\ref{eqqq:8})(\ref{eqqq:7})}, that allows us to restore, in a precise manner, the matrix~$\boldsymbol{F}$ given in~{\rm(\ref{eqqq:2})}, we have
\begin{equation}
\label{eqqq:9}
\tau\ \  \leq\ \  \sum_{i=1}^{\tau}\mathfrak{q}(T^i)\ \ \leq\ \
\begin{cases}
\tau\cdot t\; , & \text{if\/ $t$ odd}\; ,\\
\tau\cdot(t-1)\; , & \text{if\/ $t$ even}\; .
\end{cases}
\end{equation}
\end{theorem}
\noindent{}Although the exact lower and upper bounds given in~(\ref{eqqq:9}) seam to be the encouraging and discouraging ones, respectively, it would be interesting to imagine how a specific practice of preprocessing an ML-data might look like; see~Appendix.

\newpage

\section*{Appendix:\\ A precise linear algebraic portrait of a small text file}

Consider the $(\tau:=8)\times (t:=36)$ binary $\{0,1\}$-matrix $\widetilde{\boldsymbol{F}}:=$
{\tiny
\begin{equation}
\label{eqqq:10}
\left(
\begin{smallmatrix}
\phantom{-}0& \phantom{-}0& \phantom{-}0& \phantom{-}0& \phantom{-}0& \phantom{-}0& \phantom{-}0& \phantom{-}0& \phantom{-}0& \phantom{-}0& \phantom{-}0& \phantom{-}0& \phantom{-}0& \phantom{-}0& \phantom{-}0& \phantom{-}0& \phantom{-}0& \phantom{-}0& \phantom{-}0& \phantom{-}0& \phantom{-}0& \phantom{-}0& \phantom{-}0& \phantom{-}0& \phantom{-}0& \phantom{-}0& \phantom{-}0& \phantom{-}0& \phantom{-}0& \phantom{-}0& \phantom{-}0& \phantom{-}0& \phantom{-}0& \phantom{-}0& \phantom{-}0& \phantom{-}0\\
\phantom{-}1& \phantom{-}1& \phantom{-}1& \phantom{-}1& \phantom{-}0& \phantom{-}1& \phantom{-}1& \phantom{-}1& \phantom{-}0& \phantom{-}1& \phantom{-}1& \phantom{-}1& \phantom{-}1& \phantom{-}0& \phantom{-}1& \phantom{-}0& \phantom{-}1& \phantom{-}1& \phantom{-}0& \phantom{-}1& \phantom{-}1& \phantom{-}1& \phantom{-}1& \phantom{-}1& \phantom{-}0& \phantom{-}0& \phantom{-}1& \phantom{-}1& \phantom{-}1& \phantom{-}1& \phantom{-}1& \phantom{-}1& \phantom{-}1& \phantom{-}1& \phantom{-}1& \phantom{-}0\\
\phantom{-}0& \phantom{-}1& \phantom{-}1& \phantom{-}1& \phantom{-}1& \phantom{-}1& \phantom{-}1& \phantom{-}1& \phantom{-}1& \phantom{-}1& \phantom{-}1& \phantom{-}1& \phantom{-}1& \phantom{-}1& \phantom{-}1& \phantom{-}1& \phantom{-}1& \phantom{-}1& \phantom{-}1& \phantom{-}1& \phantom{-}1& \phantom{-}1& \phantom{-}1& \phantom{-}1& \phantom{-}1& \phantom{-}1& \phantom{-}0& \phantom{-}1& \phantom{-}1& \phantom{-}1& \phantom{-}1& \phantom{-}1& \phantom{-}1& \phantom{-}1& \phantom{-}1& \phantom{-}1\\
\phantom{-}0& \phantom{-}0& \phantom{-}1& \phantom{-}0& \phantom{-}0& \phantom{-}1& \phantom{-}0& \phantom{-}1& \phantom{-}0& \phantom{-}1& \phantom{-}1& \phantom{-}0& \phantom{-}1& \phantom{-}0& \phantom{-}0& \phantom{-}0& \phantom{-}1& \phantom{-}0& \phantom{-}0& \phantom{-}0& \phantom{-}0& \phantom{-}0& \phantom{-}0& \phantom{-}1& \phantom{-}0& \phantom{-}0& \phantom{-}0& \phantom{-}0& \phantom{-}1& \phantom{-}0& \phantom{-}0& \phantom{-}0& \phantom{-}0& \phantom{-}0& \phantom{-}0& \phantom{-}1\\
\phantom{-}1& \phantom{-}0& \phantom{-}0& \phantom{-}0& \phantom{-}0& \phantom{-}1& \phantom{-}1& \phantom{-}0& \phantom{-}0& \phantom{-}0& \phantom{-}0& \phantom{-}0& \phantom{-}1& \phantom{-}0& \phantom{-}0& \phantom{-}0& \phantom{-}0& \phantom{-}1& \phantom{-}1& \phantom{-}1& \phantom{-}1& \phantom{-}0& \phantom{-}1& \phantom{-}0& \phantom{-}1& \phantom{-}0& \phantom{-}0& \phantom{-}0& \phantom{-}0& \phantom{-}0& \phantom{-}0& \phantom{-}1& \phantom{-}1& \phantom{-}1& \phantom{-}0& \phantom{-}1\\
\phantom{-}0& \phantom{-}0& \phantom{-}1& \phantom{-}1& \phantom{-}0& \phantom{-}0& \phantom{-}1& \phantom{-}1& \phantom{-}0& \phantom{-}0& \phantom{-}0& \phantom{-}0& \phantom{-}0& \phantom{-}1& \phantom{-}1& \phantom{-}0& \phantom{-}1& \phantom{-}1& \phantom{-}1& \phantom{-}1& \phantom{-}0& \phantom{-}1& \phantom{-}0& \phantom{-}1& \phantom{-}1& \phantom{-}0& \phantom{-}1& \phantom{-}1& \phantom{-}0& \phantom{-}1& \phantom{-}1& \phantom{-}1& \phantom{-}1& \phantom{-}1& \phantom{-}0& \phantom{-}1\\
\phantom{-}0& \phantom{-}0& \phantom{-}1& \phantom{-}0& \phantom{-}0& \phantom{-}0& \phantom{-}1& \phantom{-}0& \phantom{-}0& \phantom{-}0& \phantom{-}1& \phantom{-}0& \phantom{-}0& \phantom{-}1& \phantom{-}0& \phantom{-}0& \phantom{-}0& \phantom{-}1& \phantom{-}0& \phantom{-}1& \phantom{-}0& \phantom{-}1& \phantom{-}0& \phantom{-}0& \phantom{-}0& \phantom{-}0& \phantom{-}0& \phantom{-}0& \phantom{-}1& \phantom{-}0& \phantom{-}0& \phantom{-}0& \phantom{-}1& \phantom{-}1& \phantom{-}0& \phantom{-}1\\
\phantom{-}0& \phantom{-}1& \phantom{-}0& \phantom{-}1& \phantom{-}0& \phantom{-}1& \phantom{-}1& \phantom{-}1& \phantom{-}0& \phantom{-}0& \phantom{-}0& \phantom{-}1& \phantom{-}1& \phantom{-}1& \phantom{-}0& \phantom{-}0& \phantom{-}0& \phantom{-}1& \phantom{-}1& \phantom{-}0& \phantom{-}1& \phantom{-}1& \phantom{-}0& \phantom{-}0& \phantom{-}0& \phantom{-}0& \phantom{-}0& \phantom{-}1& \phantom{-}1& \phantom{-}0& \phantom{-}1& \phantom{-}1& \phantom{-}1& \phantom{-}0& \phantom{-}1& \phantom{-}1
\end{smallmatrix}
\right)
\end{equation}
}
\noindent{}If we take the view of the matrix~$\widetilde{\boldsymbol{F}}$, where its columns are interpreted as `transposed' {\em bytes\/} representing (in {\em ASCII/UTF-8} character encoding) letters of the English alphabet and other things, then the binary array~(\ref{eqqq:10}) hides in itself the sharp question:
{\tiny
\begin{equation}
\label{eqqq:14}
\phantom{\big(}
\begin{smallmatrix}
\phantom{-}\phantom{0}& \phantom{-}\phantom{0}& \phantom{-}\phantom{0}& \phantom{-}\phantom{0}& \phantom{-}\phantom{0}& \phantom{-}\phantom{0}& \phantom{-}\phantom{0}& \phantom{-}\phantom{0}& \phantom{-}\phantom{0}& \phantom{-}\phantom{0}& \phantom{-}\phantom{0}& \phantom{-}\phantom{0}& \phantom{-}\phantom{0}& \phantom{-}\phantom{0}& \phantom{-}\phantom{0}& \phantom{-}\phantom{0}& \phantom{-}\phantom{0}& \phantom{-}\phantom{0}& \phantom{-}\phantom{0}& \phantom{-}\phantom{0}& \phantom{-}\phantom{0}& \phantom{-}\phantom{0}& \phantom{-}\phantom{0}& \phantom{-}\phantom{0}& \phantom{-}\phantom{0}& \phantom{-}\phantom{0}& \phantom{-}\phantom{0}& \phantom{-}\phantom{0}& \phantom{-}\phantom{0}& \phantom{-}\phantom{0}& \phantom{-}\phantom{0}& \phantom{-}\phantom{0}& \phantom{-}\phantom{0}& \phantom{-}\phantom{0}& \phantom{-}\phantom{0}& \phantom{-}\phantom{0}
\\
\phantom{-}\text{H}& \phantom{-}\text{a}& \phantom{-}\text{v}& \phantom{-}\text{e}& \phantom{-} & \phantom{-}\text{y}& \phantom{-}\text{o}& \phantom{-}\text{u}& \phantom{-} & \phantom{-}\text{p}& \phantom{-}\text{r}& \phantom{-}\text{a}& \phantom{-}\text{y}& \phantom{-}\text{'}& \phantom{-}\text{d}& \phantom{-} & \phantom{-}\text{t}& \phantom{-}\text{o}& \phantom{-}\text{-}& \phantom{-}\text{n}& \phantom{-}\text{i}& \phantom{-}\text{g}& \phantom{-}\text{h}& \phantom{-}\text{t}& \phantom{-}\text{,}& \phantom{-} & \phantom{-}\text{D}& \phantom{-}\text{e}& \phantom{-}\text{s}& \phantom{-}\text{d}& \phantom{-}\text{e}& \phantom{-}\text{m}& \phantom{-}\text{o}& \phantom{-}\text{n}& \phantom{-}\text{a}& \phantom{-}\text{?}
\end{smallmatrix}
\phantom{\big)}
\end{equation}
}

Let us associate with the matrix/array $\widetilde{\boldsymbol{F}}$, by means of the conversion
\begin{equation*}
\{0,1\}\to\{1,-1\}: \ \ \ 0\mapsto 1\; ,\ \ \ 1\mapsto -1\; ,
\end{equation*}
the $\tau\times t$ matrix/array $\boldsymbol{F}:=$
{\tiny
\begin{equation}
\label{eqqq:11}
\left(
\begin{smallmatrix}
\phantom{-}1& \phantom{-}1& \phantom{-}1& \phantom{-}1& \phantom{-}1& \phantom{-}1& \phantom{-}1& \phantom{-}1& \phantom{-}1& \phantom{-}1& \phantom{-}1& \phantom{-}1& \phantom{-}1& \phantom{-}1& \phantom{-}1& \phantom{-}1& \phantom{-}1& \phantom{-}1& \phantom{-}1& \phantom{-}1& \phantom{-}1& \phantom{-}1& \phantom{-}1& \phantom{-}1& \phantom{-}1& \phantom{-}1& \phantom{-}1& \phantom{-}1& \phantom{-}1& \phantom{-}1& \phantom{-}1& \phantom{-}1& \phantom{-}1& \phantom{-}1& \phantom{-}1& \phantom{-}1\\
-1& -1& -1& -1& \phantom{-}1& -1& -1& -1& \phantom{-}1& -1& -1& -1& -1& \phantom{-}1& -1& \phantom{-}1& -1& -1& \phantom{-}1& -1& -1& -1& -1& -1& \phantom{-}1& \phantom{-}1& -1& -1& -1& -1& -1& -1& -1& -1& -1& \phantom{-}1\\
\phantom{-}1& -1& -1& -1& -1& -1& -1& -1& -1& -1& -1& -1& -1& -1& -1& -1& -1& -1& -1& -1& -1& -1& -1& -1& -1& -1& \phantom{-}1& -1& -1& -1& -1& -1& -1& -1& -1& -1\\
\phantom{-}1& \phantom{-}1& -1& \phantom{-}1& \phantom{-}1& -1& \phantom{-}1& -1& \phantom{-}1& -1& -1& \phantom{-}1& -1& \phantom{-}1& \phantom{-}1& \phantom{-}1& -1& \phantom{-}1& \phantom{-}1& \phantom{-}1& \phantom{-}1& \phantom{-}1& \phantom{-}1& -1& \phantom{-}1& \phantom{-}1& \phantom{-}1& \phantom{-}1& -1& \phantom{-}1& \phantom{-}1& \phantom{-}1& \phantom{-}1& \phantom{-}1& \phantom{-}1& -1\\
-1& \phantom{-}1& \phantom{-}1& \phantom{-}1& \phantom{-}1& -1& -1& \phantom{-}1& \phantom{-}1& \phantom{-}1& \phantom{-}1& \phantom{-}1& -1& \phantom{-}1& \phantom{-}1& \phantom{-}1& \phantom{-}1& -1& -1& -1& -1& \phantom{-}1& -1& \phantom{-}1& -1& \phantom{-}1& \phantom{-}1& \phantom{-}1& \phantom{-}1& \phantom{-}1& \phantom{-}1& -1& -1& -1& \phantom{-}1& -1\\
\phantom{-}1& \phantom{-}1& -1& -1& \phantom{-}1& \phantom{-}1& -1& -1& \phantom{-}1& \phantom{-}1& \phantom{-}1& \phantom{-}1& \phantom{-}1& -1& -1& \phantom{-}1& -1& -1& -1& -1& \phantom{-}1& -1& \phantom{-}1& -1& -1& \phantom{-}1& -1& -1& \phantom{-}1& -1& -1& -1& -1& -1& \phantom{-}1& -1\\
\phantom{-}1& \phantom{-}1& -1& \phantom{-}1& \phantom{-}1& \phantom{-}1& -1& \phantom{-}1& \phantom{-}1& \phantom{-}1& -1& \phantom{-}1& \phantom{-}1& -1& \phantom{-}1& \phantom{-}1& \phantom{-}1& -1& \phantom{-}1& -1& \phantom{-}1& -1& \phantom{-}1& \phantom{-}1& \phantom{-}1& \phantom{-}1& \phantom{-}1& \phantom{-}1& -1& \phantom{-}1& \phantom{-}1& \phantom{-}1& -1& -1& \phantom{-}1& -1\\
\phantom{-}1& -1& \phantom{-}1& -1& \phantom{-}1& -1& -1& -1& \phantom{-}1& \phantom{-}1& \phantom{-}1& -1& -1& -1& \phantom{-}1& \phantom{-}1& \phantom{-}1& -1& -1& \phantom{-}1& -1& -1& \phantom{-}1& \phantom{-}1& \phantom{-}1& \phantom{-}1& \phantom{-}1& -1& -1& \phantom{-}1& -1& -1& -1& \phantom{-}1& -1& -1
\end{smallmatrix}
\right)
\end{equation}
}
\noindent{}The rows of $\boldsymbol{F}$ are vertices of the hypercube graph~$\boldsymbol{H}((t:=36),2)$ whose distinguished symmetric cycle~$\boldsymbol{R}$ is defined by~(\ref{eqqq:5})(\ref{eqqq:6}). Keeping in mind the description of the matrix~$\boldsymbol{F}$ given in~(\ref{eqqq:2}), let us apply Proposition~5.9 from~\cite{M-JSP} (\cite[Prop.~2.4]{M-02601}) to the array~(\ref{eqqq:11}), in order to create its linear algebraic portrait~(\ref{eqqq:12}) specified in Theorem~\ref{proppp:1}:
\begin{multline}
\label{eqqq:15}
\quad\quad\quad\quad\quad\quad\quad\quad\quad\quad\quad\quad\quad\Big(\ 36\; ,\ 8\; ,\\
\{0\}\; ,\\
\{4,8,13,15,18,24,35,41,45,50,52,55,62\}\; ,\\
\{26,37,63\}\; ,\\
\{3,6,8,11,13,17,24,29,38,41,43,45,48,52,59,64,71\}\; ,\\
\{1,7,13,21,23,25,34,36,41,48,53,58,60,67,71\}\; ,\\
\{4,8,15,20,22,25,28,34,38,42,49,52,57,59,62,65,71\}\; ,\\
\{3,7,11,14,18,20,22,29,34,38,42,46,49,53,55,57,64,68,71\}\; ,\\
\{2,4,8,14,19,22,29,33,37,39,41,47,53,56,63,66,70\}
\ \Big)\; .
\end{multline}
For example,
\begin{equation*}
\{k_0(T^3),k_1(T^3),k_2(T^3)\}=\{26,37,63\}\; ,
\end{equation*}
\newpage
and since we have $R^{26}=$
{\tiny
\begin{equation*}
\left(
\begin{smallmatrix}
-1& -1& -1& -1& -1& -1& -1& -1& -1& -1& -1& -1& -1& -1& -1& -1& -1& -1& -1& -1& -1& -1& -1& -1& -1& -1& \phantom{-}1& \phantom{-}1& \phantom{-}1& \phantom{-}1& \phantom{-}1& \phantom{-}1& \phantom{-}1& \phantom{-}1& \phantom{-}1& \phantom{-}1
\end{smallmatrix}
\right)\; ,
\end{equation*}
}
and $R^{37}=$
{\tiny
\begin{equation*}
\left(
\begin{smallmatrix}
\phantom{-}1& -1& -1& -1& -1& -1& -1& -1& -1& -1& -1& -1& -1& -1& -1& -1& -1& -1& -1& -1& -1& -1& -1& -1& -1& -1& -1& -1& -1& -1& -1& -1& -1& -1& -1& -1
\end{smallmatrix}
\right)\; ,
\end{equation*}
}
and $R^{63}=$
{\tiny
\begin{equation*}
\left(
\begin{smallmatrix}
\phantom{-}1& \phantom{-}1& \phantom{-}1& \phantom{-}1& \phantom{-}1& \phantom{-}1& \phantom{-}1& \phantom{-}1& \phantom{-}1& \phantom{-}1& \phantom{-}1& \phantom{-}1& \phantom{-}1& \phantom{-}1& \phantom{-}1& \phantom{-}1& \phantom{-}1& \phantom{-}1& \phantom{-}1& \phantom{-}1& \phantom{-}1& \phantom{-}1& \phantom{-}1& \phantom{-}1& \phantom{-}1& \phantom{-}1& \phantom{-}1& -1& -1& -1& -1& -1& -1& -1& -1& -1
\end{smallmatrix}
\right)\; ,
\end{equation*}
}
indeed, we see that $R^{26}+R^{37}+R^{63}=T^3=$
{\tiny
\begin{equation*}
\left(
\begin{smallmatrix}
\phantom{-}1& -1& -1& -1& -1& -1& -1& -1& -1& -1& -1& -1& -1& -1& -1& -1& -1& -1& -1& -1& -1& -1& -1& -1& -1& -1& \phantom{-}1& -1& -1& -1& -1& -1& -1& -1& -1& -1
\end{smallmatrix}
\right)\; .
\end{equation*}
}

For the portrait~(\ref{eqqq:15}) we have $\sum_{i=1}^{8}\mathfrak{q}(T^i)=102$.

We end this Appendix by transposing the columns of the matrices~(\ref{eqqq:10}) and~(\ref{eqqq:11}), and by gluing them into two row vectors,
$\widetilde{L}\in\{0,1\}^{\tau\cdot t}\subset\mathbb{R}^{\tau\cdot t}$ and
\begin{equation*}
L:=\left(
\begin{smallmatrix}
F^1
& F^2
& \cdots
& F^t
\end{smallmatrix}
\right)
\in\{1,-1\}^{\tau\cdot t}\subset\mathbb{R}^{\tau\cdot t}
\; ,
\end{equation*}
respectively, of dimension~$\tau\cdot t:=8\cdot 36=288$; cf. description~(\ref{eqqq:2}) of the matrix~$\boldsymbol{F}$. Thus,~$\widetilde{L}:=$
{\tiny
\begin{equation}
\label{eqqq:13}
\left(
\begin{smallmatrix}
\phantom{-}0& \phantom{-}1& \phantom{-}0& \phantom{-}0& \phantom{-}1& \phantom{-}0& \phantom{-}0& \phantom{-}0&
\phantom{-}0& \phantom{-}1& \phantom{-}1& \phantom{-}0& \phantom{-}0& \phantom{-}0& \phantom{-}0& \phantom{-}1& & &
\cdots& \cdots& \cdots&
\phantom{-}0& \phantom{-}1& \phantom{-}1& \phantom{-}0& \phantom{-}0& \phantom{-}0& \phantom{-}0& \phantom{-}1&
\phantom{-}0& \phantom{-}0& \phantom{-}1& \phantom{-}1& \phantom{-}1& \phantom{-}1& \phantom{-}1& \phantom{-}1
\end{smallmatrix}
\right)\; ,
\end{equation}
}
and $L:=$
{\tiny
\begin{equation}
\label{eqqq:16}
\left(
\begin{smallmatrix}
\phantom{-}1& -1& \phantom{-}1& \phantom{-}1& -1& \phantom{-}1& \phantom{-}1& \phantom{-}1&
\phantom{-}1& -1& -1& \phantom{-}1& \phantom{-}1& \phantom{-}1& \phantom{-}1& -1& & &
\cdots& \cdots& \cdots&
\phantom{-}1& -1& -1& \phantom{-}1& \phantom{-}1& \phantom{-}1& \phantom{-}1& -1&
\phantom{-}1& \phantom{-}1& -1& -1& -1& -1& -1& -1
\end{smallmatrix}
\right)\; ,
\end{equation}
}
Now, the matrix~$\widetilde{L}$ given in~(\ref{eqqq:13}) is essentially a sequence of {\em bits\/} composing a text {\em file\/} that consists of the {\em string\/} (i.e., a sequence of {\em characters}) displayed in~(\ref{eqqq:14}). By applying Proposition~5.9 from~\cite{M-JSP} (\cite[Prop.~2.4]{M-02601}) to the vertex~(\ref{eqqq:16}) of the hypercube graph~$\boldsymbol{H}(288,2)$, we obtain the following exact linear algebraic portrait of the string (in fact, of the file containing the string) displayed in~(\ref{eqqq:14}):
\begin{multline}
\label{eqqq:17}
\quad\quad\quad\quad\quad\quad\quad\quad\quad\quad\quad\quad\quad\Big(\ 288\; ,\ 1\; ,\\
\{2,5,11,16,20,23,27,30,32,35,45,48,51,56,60,62,64,67,76,84,87,91,\\
96,101,104,107,112,115,118,123,132,134,139,144,147,150,152,155,\\
159,163,165,168,171,176,179,181,188,190,195,198,203,210,214,219,\\
222,224,228,232,235,238,243,246,248,251,254,256,259,264,267,271,\\
275,280,289,292,297,303,305,309,313,317,319,322,329,335,337,340,\\
345,349,351,354,361,369,374,377,383,385,391,394,397,401,405,410,\\
417,421,425,428,434,436,439,441,444,449,452,455,457,461,465,468,\\
473,477,482,484,490,497,501,505,509,511,513,518,521,525,529,533,\\
535,537,540,543,545,548,553,556,561,567,570\}
\ \Big)\; .
\end{multline}
For the portrait~(\ref{eqqq:17}) we have $\mathfrak{q}(L):=|\boldsymbol{Q}(L,\boldsymbol{R})|=145$.

\vspace{5mm}
\end{document}